\documentclass[12pt]{article}
\usepackage{geometry}
\usepackage{graphicx}

\usepackage{url}
\begin{document}

\begin{center}

Author's accepted manuscript. Please cite: Isobel Falconer, `Historical Notes: The Gravitational Constant', {\em Mathematics Today}, vol.58 no.4 (August 2022) pp126-127
\vspace{5mm} %5mm vertical space

\bf The Gravitational Constant

\bf Isobel Falconer
\end{center}

{\em It is probably within the knowledge of most of those present that Sir Isaac Newton, by his great discovery of gravitation and its laws, was able to show that a single principle … viz. that every particle in the universe attracts any other particle towards itself with a force which is proportional to the product of their masses divided by the square of the distance between them}, declaimed Charles Vernon Boys in 1894 \cite{Boys}.

Boys had invented a method of extruding fine quartz threads, up to 60 feet long, by using a crossbow to shoot a straw arrow from a stock of molten quartz. He used the threads to suspend the torsion balance with which he determined the gravitational constant, ‘$G$’ with an unprecedented precision of five significant figures.     

‘$G$’ appears in ‘Newton’s law of gravitation’ as we usually write it today: 
\begin{equation}
F = G\frac{Mm}{r^2}
\end{equation}
where $F$ is the force due to gravity between two masses $M$, $m$ a distance $r$ from each other. Or, if $m$ is the earth, $\Delta$ its mean density, and $M$ a mass near the surface, (1) can be re-written as $$F=G\frac{4\pi}{3}rM\Delta$$ According to Wikipedia – a useful barometer of popular belief – $G$ is a ‘notable example’ of a fundamental physical constant, and among the ‘most widely recognized’ \cite{Wiki}.  

I was surprised, then, to find that over 200 years after Newton published the {\em Principia} \cite{Newt}, Boys felt the need to argue strongly for $G$: {\em G… represents that mighty principle under the influence of which every star, planet and satellite in the universe pursues its allotted course… It is in no way dependent on the accidental size or shape of the earth; if the solar system ceased to exist it would remain unchanged…} \cite{Boys}. 

Why was Boys so vehement? Six years later, in John Henry Poynting’s 1900 lecture ‘Recent Studies in Gravitation’ we find out: {\em Professor Boys has almost indignantly disclaimed that he was engaged in any such purely local experiment as the determination of the mean density of the earth. He was working for the Universe, seeking the value of G, information which would be as useful on Mars or Jupiter or out in the stellar system as here on the earth. But perhaps we may this evening consent to be more parochial in our ideas, and express the results in terms of the mean density of the earth} \cite{Poyn}.

This is the clue. $G$ was a very recent innovation; until 1884 finding the mean density of the earth was the main purpose given for gravitational experiments, dominating gravitational writing throughout the century as suggested by the Ngrams in Figure 1. Indeed Poynting, Professor of Physics at Birmingham, acquired popular fame as ‘the man who weighed the earth’ - although we remember him these days for the ‘Poynting vector’ describing the flow of energy in an electromagnetic field.

So, $G$ did not appear until 200 years after Newton proposed his gravitational theory.  Newton himself had not used symbolic equations at all, for example \cite[Book 1 Prop 69 Thm 29]{Newt}: {\em absolute forces of the attracting bodies A and B will be to each other as those very bodies A and B to which those forces belong} and \cite[Book3 Prop 5 Thm 5 Cor II]{Newt}{\em The force of gravity which tends to any one planet is inversely as the square of the distance of places from that planet’s centre} \. By comparing the attraction of one planet to that of another, expressing the results as ratios, there was no place for a gravitational constant. In this relative system, finding the mass or mean density of the earth was a necessary first step to establishing the masses of astronomical bodies.

Around 1800 two factors began to change this situation. First, in 1798, the reclusive chemist, Henry Cavendish, realised the first laboratory measurement of the attraction between known masses. All the while gravitational calculations depended on astronomical observations it would have been impossible to measure $G$ and $M$ independently of each other even if the equations had been cast in our modern form; only the product $GM$ could be measured and the question of a constant $G$ did not arise.

Cavendish used a torsion balance to compare the attraction of a lead weight to a lead ball with the attraction of the ball to the earth. Figure 2 shows his apparatus, housed in a special building and observed remotely from outside using a lamp and telescope. Using relative measures, he concluded that {\em the attraction of the leaden weight on the ball will be to that of the earth thereon, as $10,64 \, \times \, ,9779 \times (6/8,85)^2$ to $41 800 000D$}, where $10.64 \times 0.9779 \times (6/8.85)^2$ is his observational result for the force of attraction between weight and ball, 41800000 is the diameter of the earth in feet, and $D$ is the density of the earth relative to that of water \cite{Cave}. Cavendish’s use of a comma instead of a decimal point shows that point notation was not yet standard in Britain.  The ‘diameter of the earth’ is implicitly relative to that of a 1 foot sphere of water so that both sides of the ratio are dimensionless as was conventional at the time.

But even now the question of $G$ would not have arisen without Euler’s introduction of functional notation in 1734. In functional notation, the dependent variable – in this case the gravitational force on a body – was typically written on the left hand side of an equation – with the requisite combination of independent variables on the right. This form lent itself to absolute, rather than relative, measurements - and crucially a coefficient of proportionality was required to satisfy requirements of dimensionality and ‘make the units come out right’. For gravitation we see this in Laplace’s {\em Mecanique Celeste} (1798):
$$ \phi = \frac{4\pi^2}{k^2}\frac{1}{r^2}$$ where $\phi$ was the accelerating force between two bodies and $k$ was the proportionality constant\cite{Lapl}.
				
It took a while for the British to catch on, but the Cambridge mathematician John Pratt did so in his 1836 textbook designed to introduce students to Continental mathematics: {\em Let $M$ and $m$ be the masses of the two particles, $r$ their distance ... then, if the unit of attraction be the attraction of a unit of mass at a unit of distance, the accelerating force produced in $M$ by the attraction of $m=m/r^2$} \cite{Prat}. No constant appeared explicitly, as Pratt had carefully defined his units so that the coefficient of proportionality was 1. It looks as though his aim was to convey clearly the functional mathematical form without confusing the issue with physical realism.

Real life experimenters were denied this luxury. Frances Baily and George Biddell Airy, the Astronomer Royal, had to introduce a coefficient when using a torsion balance to measure the density of the earth in 1843 \cite[pp103, 109]{Airy}: {\em The unit which we shall adopt is, the mass which at London weighs one grain: and for the attraction which it produces at distance 1 … we shall multiply it by the modulus $k$}.
$$k=\frac{2\pi^2Fh}{Ei}\frac{(s-s')}{T^2}$$
where $F$, $h$, $E$, and $i$ are known instrumental parameters, $s$, $s'$ are the observed equilibrium points of the balance with the deflecting masses in two different positions, and $T$ is the observed period of oscillation.

I find Airy’s use of a ‘modulus $k$’ illuminating. A ‘modulus’ expresses proportionality but implies no commitment to constancy. For example, ‘Young’s modulus’, $E$, shows that the relationship between the elastic deformation of a wire and a small applied straining force is linear, but $E$ varies from one wire to the next. Airy’s choice of ‘$k$’ for the modulus can be traced back at least to Legendre’s 1792 memoir on elliptic transcendentals \cite{Lege} – evidence for its origin in mathematical questions that ultimately derived from astronomical calculations.
 
In 1873, Alfred Cornu and Jean-Baptistin Baille were the first to mention a 'constant of attraction' in \cite{Corn}, coupling it with the mean density of the Earth. In 1884, 11 years later, came the first paper that emphatically foregrounded the gravitational constant, now finally denoted by $G$, as its aim: Franz Richarz and Arthur König’s `Eine neue Methode zur Bestimmung der Gravitationsconstante' \cite{Rich}. 

Which brings us to the Poynting – Boys debate. Both men wrote the gravitational equations in our familiar form: $F=GMm/r^2$, and evidence was piling up from successive measurements of the mean density of the earth that $G$ was always the same. So was Poynting merely being ultra-conservative in rejecting Boys’ reification of $G$ as a fundamental physical constant?  On the contrary. He believed that such a commitment was an unproductive dead end - variation in $G$ should be actively sought. While Boys glorified a constant $G$, as a ‘mighty principle’ governing the universe, Poynting pointed out that constancy made gravity unlike any other physical phenomenon; reasoning by analogy from phenomena such as electrostatics whose laws had a similar mathematical form would be impossible: {\em this unlikeness, this independence of gravitation of any quality but mass, bars the way to any explanation of its nature… Gravitation still stands alone} \cite{Poyn}. Despite general relativity, gravitation still stands alone today.

\newpage

{\bf Figures}

\vspace{10mm} %10mm vertical space

Figure 1: A Google Ngram comparison for gravity phrases 1800-1900

\includegraphics[scale=0.5]{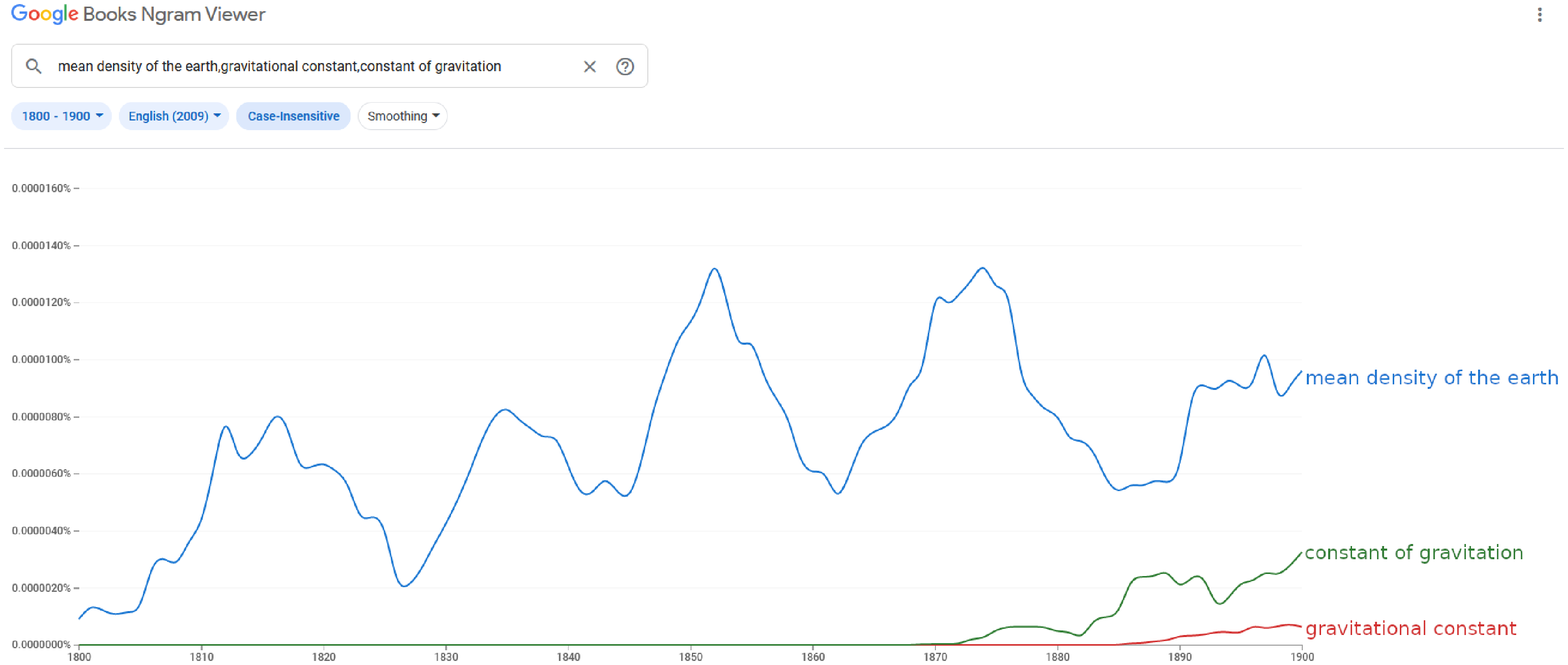}

\vspace{10mm} %10mm vertical space

Figure 2: Cavendish's diagram of his torsion apparatus for measuring the mean density of the earth \cite{Cave}

\includegraphics{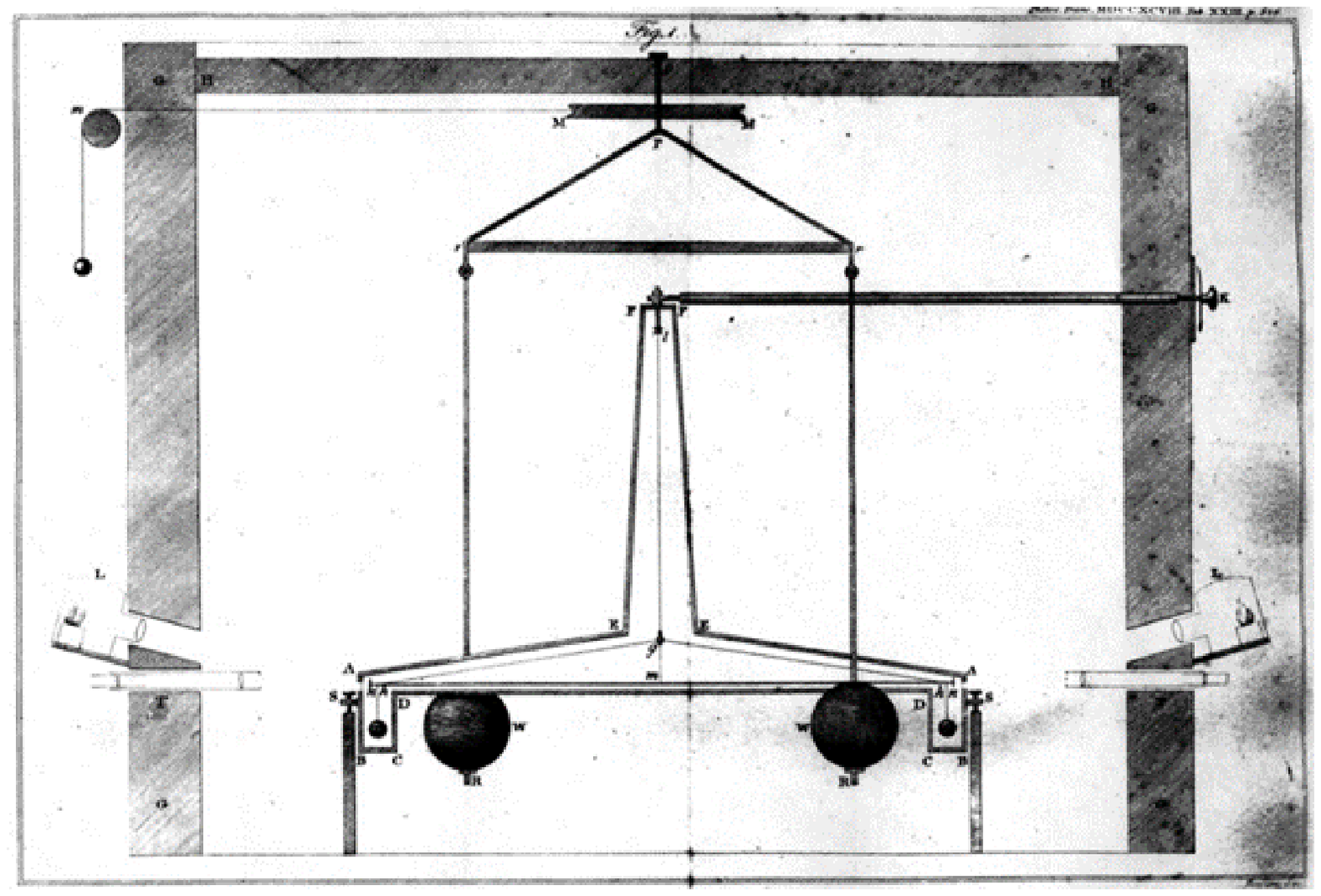}


\begin{thebibliography}{11}

\bibitem{Boys}
C.V. Boys, `The Newtonian Constant of Gravitation', Address to the Royal Institution of Great Britain, (8 Jun. 1894) 

\bibitem{Wiki}
Wikipedia, `Physical constant' \url{en.wikipedia.org/wiki/Physical_constant}

\bibitem{Poyn}
J.H. Poynting, `Recent Studies in Gravitation', Address to the Royal Institution of Great Britain (23 Feb. 1900)

\bibitem{Newt}
I. Newton, {\em Mathematical Principles of Natural Philosophy }, trans. A. Motte (1729)

\bibitem{Cave}
H. Cavendish, `Experiments to determine the density of the earth' {\em Phil. Trans. R. Soc.} 88 (1798) 469–526 

\bibitem{Lapl}
P.S. Laplace, {\em Mécanique Céleste} trans. N. Bowditch (1829)

\bibitem{Prat}
H. Pratt, {\em The Mathematical Principles of Mechanical Philosophy} (1836)

\bibitem{Airy}
G.B. Airy, `On the mathematical theory of Cavendish's Experiment', {\em Memoirs of the Royal Astronomical Society} 14 (1843) 99-110

\bibitem{Lege}
A-M. Legendre, {\em Mémoire sur les transcendantes elliptiques} (1792)

\bibitem{Corn}
M.A. Cornu and J.-B. Baille, `Détermination nouvelle de la constante de l’attraction et de la densité moyenne de la terre', {\em Comptes Rendus de l'Académie des sciences} 76 (1873) 954-957

\bibitem{Rich}
A. König and F. Richarz, `Eine neue Methode zur bestimmung der Gravitationsconstante', {\em Annalen der Physik} 260 (1885) 664-668

\end{thebibliography}
\end{document}